\documentclass[12pt]{article}

\usepackage{amsmath,amssymb,amsthm,amsfonts}
\usepackage{nicematrix,tikz}
\usetikzlibrary{fit}
\tikzset{highlight/.style={rectangle,
fill=gray!50,
rounded corners = 0.5 mm, 
inner sep=1pt,
fit=#1}}

\def\Hom{{\rm Hom}}
\def\ind{{\rm ind}}
\def\ad{\mathop{\rm ad}}
\def\im{\mathop{\rm im}}

\def\phi{\varphi}
\def\g{\mathfrak g}

\def\h{\mathfrak h}
\def\F{\mathbb F}

\newtheorem{theorem}{Theorem}[section]
\newtheorem{lemma}[theorem]{Lemma}

\newtheorem{definition}[theorem]{Definition}

\theoremstyle{remark}

\providecommand{\keywords}[1]{\noindent{Keywords:} #1}
\providecommand{\classify}[1]{\noindent{Mathematics Subject
    Classification:} #1}

\title{Restricted Central Extensions of Restricted Heisenberg Lie Superalgebras}

\author{
Yong Yang\\
College of Mathematics and System Science\\
Xinjiang University\\
Urumqi 830046, China\\
yangyong195888221@163.com
}

\date{}

\begin{document}
\maketitle

\begin{abstract}
Restricted Heisenberg Lie
  superalgebras are studied over an algebraically closed field $\F$ of
  characteristic $p>2$. We use the ordinary 1- and 2-cohomology spaces with
  trivial coefficients to compute the restricted 2-cohomology
  spaces. As an application,  the corresponding 
  restricted one-dimensional central extensions are described. 
\end{abstract}

{\footnotesize

\classify{17B50; 17B56}}

\section{Introduction}
Given a symplectic form on a symplectic space, we can define its Heisenberg Lie algebra. It is
a two-step nilpotent Lie algebra with a one-dimensional center. Due to its application in the commutative
relations of quantum mechanics, Heisenberg Lie algebras have become an important research
object in modern mathematics. The concept of Heisenberg Lie algebras was generalized to Lie superalgebras by considering its supersymplectic form,
which is called Heisenberg Lie superalgebra. That is, a two-step nilpotent Lie superalgebra with a
one-dimensional center. 
In \cite{R}, the authors showed
the following classification theorem.

\begin{theorem}(see \cite[Section 3]{R})
Finite-dimensional complex Lie superalgebras of Heisenberg type are classified into two types described as follows (we write the even basis elements first separated from the odd
ones by a vertical bar):

$(1)$  $\mathfrak{h}_{m,n}$ is the $(2m+1,n)$-dimensional vector superspace spanned by the elements \[\{x_1,\ldots , x_{2m},x_{2m+1}\mid y_1,\ldots, y_n\}\] with the non-vanishing Lie
  brackets
  \[[x_i,x_{m+i}]=[y_j,y_j]=x_{2m+1},\] for $1\le i\le m$ and $1\le j\le n$.

$(2)$  $\mathfrak{ba}_{n}$ is the $(n,n+1)$-dimensional vector superspace spanned by the elements \[\{x_1,\ldots , x_{n}\mid y_1,\ldots, y_n,y_{n+1}\}\] with the non-vanishing Lie
  brackets
  \[[x_i,y_i]=y_{n+1},\] for $1\le i\le n$.
\end{theorem}

The cohomology theory of Heisenberg Lie algebras and Heisenberg Lie superalgebras was developed by many authors, see \cite{BaLi,CJ,EFY,S,SK} for examples.
In particular, Bai and Liu   computed the trivial cohomology of Heisenberg Lie superalgebras over fields of characteristic zero in \cite{BaLi}.

Based on the fact that the $p$-powers of inner derivations of associative algebras are still inner derivations, the notion of restricted Lie algebras over ﬁelds of characteristic $p>0$ was introduced and Jacobson proved that a Lie algebra is restrictable if and only if the $p$-powers of the inner derivations are still inner derivations, see  \cite{J1,J2}. 
The notion was generalized to Lie superalgebras called restricted Lie superalgebras or Lie $p$-superalgebras in \cite{M}. Different from restricted Lie algebras, a Lie superalgebra is restrictable if and only if 
the $p$-powers of the even inner derivations are still inner derivations \cite{P}. As a result, two-step nilpotent Lie superalgebras, particularly Heisenberg Lie superalgebras, are restrictable.
The restricted cohomology theory of restricted Lie algebras was first considered by Hochschild in \cite{H}.
For a restricted Lie algebra, he associated  the usual Chevalley-Eilenberg cohomology with the restricted cohomology  by the canonical projection from its universal enveloping algebra
to the restricted universal enveloping algebra. Moreover, Evans and Fuchs constructed a cochain complex for restricted cohomology up to order 3 in \cite{EFu}. It has been proved to be an effective tool to compute restricted cohomology.
By means of the cochain complex, Evans and Fialowski, etc. computed the restricted cohomology and restricted central extensions for restricted Witt algebras, restricted affine nilpotent Lie algebras, restricted simple Lie algebras, restricted filiform Lie algebras and restricted Heisenberg Lie algebras, see \cite{EFi1,EFi2,EFi3,EFi4,EFP,EFY}. 
Moreover, the central extensions of symmetric modular Lie superalgebras
 have been studied in \cite{BG}. Recently, the restricted cohomology theory of restricted Lie superalgebras was introduced by Bouarroudj and Ehret \cite{BE}.

 In this paper, restricted Heisenberg Lie superalgebras are studied.
We compute the restricted 1-, 2-cohomology spaces with trivial coeﬃcients and the corresponding restricted one-dimensional central extensions of restricted Heisenberg Lie
superalgebras. The structure of the paper is as follows.  In section 2, we recall the definitions
 of restricted Lie superalgebras, the Chevalley-Eilenberg cochain complex and restricted cochain
complex  in the case of trivial coefficients. We also recall the six-term exact sequence associated the ordinary cohomology with restricted cohomology and the correspondence between restricted one
dimensional central extensions and restricted 2-cohomology.
In section 3, we compute the restricted cohomology and restricted one-dimensional central extensions of
$\h^{\lambda}_{m,n}$, the restricted Heisenberg Lie superalgebra with the even center.
In section 4, we compute the restricted cohomology and restricted one-dimensional central extensions of
$\mathfrak{ba}^{0}_n$, the restricted Heisenberg Lie superalgebra with the odd center.
We only consider central extensions obtained by even cocycles, because the central extensions obtained by the odd ones are trivial.

\section{Restricted Lie superalgebras and cohomology}

In this section, we recall the definitions of a restricted Lie superalgebra,
the corresponding Chevalley-Eilenberg cochain complex and restricted cochain complex \cite{BE}.
 Everywhere in this
section, $\F$ denotes an algebraically closed field of characteristic
$p>2$ and $\g$ denotes a finite-dimensional Lie superalgebra over $\F$ with an
ordered basis $\{e_1,\ldots,e_m\mid e_{m+1},\ldots,e_{m+n}\}$ (we
write the even basis elements first separated from the odd ones by a vertical
bar).
Write $e^k$,  $e^{i,j}$ and $e^{u,v,w}$ for the dual vectors of the basis vectors $e_k\in\g$, $e_{i,j}=e_i\wedge e_j\in \wedge^{2}\g$ and $e_{u,v,w}=e_u\wedge e_v\wedge e_w\in\wedge^{3}\g$, respectively.
In the case of $p=3$,  we add $[g,[g,g]]=0$ for all $g\in\g_{\bar{1}}$ as a part of the definition of Lie superalgebras.
For $j\geq2$ and $g_1,\ldots, g_{j}\in\g$, we denote the $j$-fold bracket
$$[g_1, g_2, g_3,\ldots,g_j] = [[\ldots[[g_1, g_2], g_3],\ldots], g_j].
$$

\subsection{Restricted Lie superalgebras}
A \emph{restricted Lie algebra}  $\mathfrak{g}$ over a field $\F$ of positive characteristic $p$ is a Lie
algebra $\mathfrak{g}$ over $\F$ together with a map
$[p]:\g\to\g$, written $g\mapsto g^{[p]}$, such that for all
$a\in\F$, and all $g,h\in \mathfrak{g}$
\begin{itemize}
\item[(1)] $(a g)^{[p]}=a^{p}g^{[p]}$;
\item[(2)]
  $(g+h)^{[p]}=g^{[p]}+h^{[p]}+ \sum\limits_{i=1}^{p-1} s_i(g,h)$,
  where $is_i(g,h)$ is the coefficient of $t^{i-1}$ in the formal expression
  $(\mathrm{ad}(tg+h))^{p-1}(g)$; and
\item[(3)] $(\mathrm{ad}\ g)^{p}=\mathrm{ad}\ g^{[p]}$.
\end{itemize}

The map $[p]$ is called a \emph{$[p]$-operator} on
$\mathfrak{g}$.  Let $(\g,[p])$ be a restricted Lie algebra. A $\g$-module $M$ is called \emph{restricted} if, for $g\in\g$ and $m\in M$,
\[\underbrace{g\ldots g}_{p}\cdot m=g^{[p]}\cdot m.\]
\begin{definition}\label{res-Lie-super}
A \emph{restricted Lie superalgebra} is a Lie superalgebra $\g=\g_{\bar{0}}\oplus \g_{\bar{1}}$ with a map $[p]:\g_{\bar{0}} \rightarrow \g_{\bar{0}}$ such that  $(\g_{\bar{0}},[p])$ is a restricted Lie algebra and $\g_{\bar{1}}$ is a restricted $\g_{\bar{0}}$-module with respect to the bracket.
\end{definition}
A Lie superalgebra $\g$ is called \emph{restrictable} if there exists a map $[p]:\g_{\bar{0}} \rightarrow \g_{\bar{0}}$ such that $(\g,[p])$ is a restricted Lie superalgebra. The superization of Jacobson’s theorem \cite{BE,J1,J2} shows that a Lie superalgebra $\g$ with  a basis $\{e_1,\ldots,e_m\mid e_{m+1},\ldots,e_{m+n} \}$
 is
restrictable if and only if $(\ad e_i)^p$ is an inner derivation for
all $1\le i\le m$. As a result,  all two-step nilpotent Lie superalgebras, particularly Heisenberg Lie superalgebras, are restrictable.


\subsection{Chevalley-Eilenberg Lie superalgebra cohomology}

We only describe  the
Chevalley-Eilenberg cochain spaces $C^q(\g)=C^q(\g,\F)$
for $q=0,1,2,3$ and differentials $d^q:C^q(\g)\to C^{q+1}(\g)$
for $q=0,1,2$.
The full description of the complex can be found in \cite{ChE,F}. Set
$C^0 (\g)=\F$ and $C^q (\g)=C^q (\g)_{\bar{0}}\oplus C^q (\g)_{\bar{1}}=
(\wedge^q\g^*)_{\bar{0}}\oplus (\wedge^q\g^*)_{\bar{1}}$ for $q=1,2,3$. 
The differentials 
$$d^q:C^q (\g)\to C^{q+1}(\g)$$ 
are defined for $\psi\in C^1 (\g)$,
$\phi\in C^2 (\g)$ and $g,h,f\in\g$ by\small
\begin{align*}
  d^0: C^0 (\g)\to C^1 (\g),  &\  d^0=0&\\
  d^1:C^1 (\g)\to C^2 (\g), &\   d^1(\psi)(g\wedge h)=\psi([g,h])&\\
  d^2:C^2 (\g)\to C^3 (\g), &\  d^2(\phi)(g\wedge h\wedge f)=\phi([g,h]\wedge f)-(-1)^{|f||h|}\phi([g,f]\wedge h)&\\
&\qquad \qquad\qquad  \qquad  \        +(-1)^{|g|(|h|+|f|)}\phi([h,f]\wedge g). &
\end{align*}
The maps $d^q$ satisfy $d^{q}d^{q-1}=0$ and 
$H^q(\g)=H^q(\g,\F)=Z^q(\g,\F)/B^q(\g,\F)$,
where $Z^q(\g,\F)=\ker(d^q)$ and $B^q(\g,\F)=\im(d^{q-1})$. A cochain is called a \emph{cocycle (resp. coboundary)} if it is in  $Z^q(\g,\F)$ (resp. $B^q(\g,\F)$).

\subsection{Restricted Lie superalgebra cohomology}
 
Since we are interested
 in central extensions, we only consider trivial coefficients and the full description can be found in
\cite{BE}.

Given $\phi\in C^2(\g)$, a map $\omega:\g_{\bar{0}}\to\F$ is \emph{
  $\phi$-compatible} if for all $g,h\in\g_{\bar{0}}$ and all $a\in\F$
\\
  
$\omega(a g)=a^p \omega (g)$ and
\begin{equation}
  \label{starprop}
  \omega(g+h)=\omega(g)+\omega(h) + \sum_{\substack{g_i=\mbox{\rm\scriptsize $g$
        or $h$}\\ g_1=g, g_2=h}}
  \frac{1}{\#(g)}\phi([g_1,g_2,g_3,\dots,g_{p-1}]\wedge g_p)
\end{equation}
where $\#(g)$ is the number of factors $g_i$ equal to $g$.

Similar to the Lie algebras case (see \cite{EFi2}),
for $\phi\in C^2(\g)$, we can assign the values of $\omega$ arbitrarily on a
basis for $\g_{\bar{0}}$ and use (1) to define
$\omega: \g_{\bar{0}} \to \F $ that is $\phi$-compatible and 
the map $\omega$ is unique because its values are completely
determined by $\omega(e_i)$ and $\phi$. In
particular, given $\phi$, we can define $\tilde\phi(e_i)=0$ for all $i$ and use
(\ref{starprop}) to determine a unique $\phi$-compatible map
$\widetilde\phi:\g_{\bar{0}}\to\F$. Note that, in general, $\widetilde\phi\ne 0$
but $\widetilde\phi (0)=0$. Moreover, If $\phi_1,\phi_2\in C^2(\g)$
and $a\in\F$, then
$\widetilde{(a\phi_1+\phi_2)} = a\widetilde\phi_1 + \widetilde\phi_2$.
\\

If $\zeta\in C^3(\g)$, then a map $\eta:\g\times \g_{\bar{0}}\to\F$ is \emph{
  $\zeta$-compatible} if for all $a\in\F$ and all $g\in\g,h,h_1,h_2\in\g_{\bar{0}}$,
$\eta(\cdot,h)$ is linear in the first coordinate,
$\eta(g,a h)=a^p\eta(g,h)$ and
\begin{align*}
  \eta(g,h_1+h_2) &=
                    \eta(g,h_1)+\eta(g,h_2)-\nonumber \\
                  & \sum_{\substack{l_1,\dots,l_p=1 {\rm or} 2\\ l_1=1,
  l_2=2}}\frac{1}{\#\{l_i=1\}}\zeta (g\wedge
  [h_{l_1},\cdots,h_{l_{p-1}}]\wedge h_{l_{p}}),
\end{align*}
where $\#\{l_i=1\}$ is the number of factors $l_i$ equal to 1. The restricted cochain spaces are defined as $C^0_{res}(\g)=C^0 (\g)$,
$C^1_{res}(\g)=C^1 (\g)$,
\[C^2_{res}(\g)=\{(\phi,\omega)\ |\ \phi\in C^2 (\g), \omega:\g_{\bar{0}}\to\F\
  \mbox{\rm is $\phi$-compatible}\},\]
\[C^3_{res}(\g)=\{(\zeta,\eta)\ |\ \zeta\in C^3 (\g),
  \eta:\g\times\g_{\bar{0}}\to\F\ \mbox{\rm is $\zeta$-compatible}\}.\] 

We define
\begin{align*}
\Hom_{\rm Fr}(\g_{\bar{0}},\F) =& \{f:\g_{\bar{0}}\to\F\ |\ f(a x+b
  y)=a^pf(x)+b^pf(y)\ \mathrm{for\ all}\ a,b\in\F\\
&\ \mathrm{and}\  x,y\in \g_{\bar{0}}\}
\end{align*}
 to be
the space of \emph{ Frobenius homomorphisms} from $\g_{\bar{0}}$ to $\F$.
For $1\le i\le m$, define $\overline e^i:\g_{\bar{0}}\to\F$ by
$\overline e^i \left(\sum_{j=1}^m a_j e_j\right ) = a_i^p.$ The set
$\{\overline e^i\ |\ 1\le i\le m\}$ is a basis for the space of Frobenius homomorphisms
$\Hom_{\rm Fr}(\g_{\bar{0}},\F)$. 

Note that a map
$\omega:\g_{\bar{0}}\to\F$ is $0$-compatible if and only if
$\omega\in \Hom_{\rm Fr}(\g_{\bar{0}},\F)$, then we have the exact sequence
\begin{align*}
  0 \longrightarrow  \Hom_{\rm Fr}(\g_{\bar{0}},\F)\stackrel{\iota}{\longrightarrow} C^2_{res}(\g) \stackrel{\pi}{\longrightarrow} C^2(\g)  \longrightarrow 0,
\end{align*}
where the maps $\iota: \Hom_{\rm Fr}(\g_{\bar{0}},\F)\to C^2_{res}(\g)$ and  $\pi: C^2_{res}(\g) \to C^2(\g)$
are given by
$$\iota(f)=(0,f)\quad \mathrm{and}\quad
\pi(\phi,\omega)=\phi,$$
for $f\in\Hom_{\rm Fr}(\g_{\bar{0}},\F)$ and $(\phi,\omega)\in C^2_{res}(\g)$. Recall $\g$ is a finite-dimensional Lie superalgebra with an
ordered basis $\{e_1,\ldots,e_m\mid e_{m+1},\ldots,e_{m+n}\}$.
Then
$\dim C^2_{res}(\g) =\dim C^2(\g)+\dim\g_{\bar{0}}$ and
\[
  \{(e^{i,j},\widetilde{e^{i,j}})\ |\ 1\le i<j \le m\} \cup \{(e^{i,j},\widetilde{e^{i,j}})\ |\ 1\le i\le m, m+1\leq j \le m+n\}\]
\[\cup
\{(e^{i,j},\widetilde{e^{i,j}})\ |\ m+1\le i\le j \le m+n\}
\cup \{(0,\overline e^i)\ |\ 1\le i\le m\}\]
is a basis for $C^2_{res}(\g)$. We will use this basis in all computations that
follow.

Define $d_{res}^0=d^0$. For $\psi\in C^1_{res}(\g)$, define the map
$\ind^1(\psi):\g_{\bar{0}}\to\F$ by
 \[\ind^1(\psi)(g)=\psi(g^{[p]}).\] 
The map
$\ind^1(\psi)$ is $d^1(\psi)$-compatible for all $\psi\in C^1_{res}(\g)$,
and the differential 
$$d^1_{res}:C^1_*(\g)\to C^2_*(\g)$$ 
is defined by
\begin{equation}
  d^1_{res}(\psi) = (d^1(\psi),\ind^1(\psi)).
\end{equation}
By definition,
\begin{equation}\label{h1}
 H^1_{res}(\g)=(\g/([\g,\g]+\langle \g_{\bar{0}}^{[p]}\rangle_\F))^*.
\end{equation}

For $(\phi,\omega)\in C^2_{res}(\g)$, define the map
$\ind^2(\phi,\omega):\g \times\g_{\bar{0}} \to\F$ by the formula
\[\ind^2(\phi,\omega)(g,h)=\phi(g\wedge h^{[p]})-\varphi([g,\underbrace{h,\ldots,h}_{p-1}],h).\]
The map $\ind^2(\phi,\omega)$ is $d^2(\phi)$-compatible for all
$\phi\in C^2(\g)$, and the differential $d^2_{res}:C^2_{res}(\g)\to C^3_{res}(\g)$
is defined by
\begin{equation}
  d^2_{res}(\phi,\omega) =
  (d^2(\phi),\ind^2(\phi,\omega)). 
\end{equation}
Note that  if $\omega_1$ and $\omega_2$ are
both $\phi$-compatible, then
$\ind^2(\phi, \omega_1)=\ind^2(\phi, \omega_2)$.
The following result is a   straightforward superization  of \cite[Lemma 2.1]{EFY}.
\begin{lemma}
  \label{swap}
  If $(\phi,\omega)\in C^2_{res}(\g)$ and $\phi=d^1(\psi)$ with
  $\psi\in C^1 (\g)$, then $(\phi,\ind^1(\psi))\in C^2_{res}(\g)$ and
  $\ind^2(\phi,\omega)=\ind^2(\phi,\ind^1(\psi))$.
\end{lemma}

These maps $d_{res}^q$ satisfy $d_{res}^{q}d_{res}^{q-1}=0$ and we define
\[H_{res}^q(\g)=H_{res}^q(\g,\F)=\ker(d_{res}^q)/\im(d_{res}^{q-1}), \]
where $q=1,2$.

\subsection{The six-term exact sequence}
Suppose that $(\g,[p])$ is a restricted Lie superalgebra. There is a six-term exact sequence, which is superized version of \cite[p. 575]{H} and \cite{Viv}:
\begin{equation}
\begin{aligned}
  \label{sixterm}
  0 &\longrightarrow H^1_{res}(\g)\stackrel{\iota_1}{\longrightarrow} H^1(\g)\stackrel{D}{\longrightarrow}\Hom_{\rm Fr}(\g_{\bar{0}},\F)   \stackrel{\iota_2}{\longrightarrow}  H^2_{res}(\g)\stackrel{\pi}{\longrightarrow}H^2(\g)\\
 & \stackrel{H}{\longrightarrow}
 \Hom_{\rm
    Fr}(\g_{\bar{0}},H^1(\g)). 
\end{aligned}
\end{equation}
The maps
$D: H^1(\g)\to\Hom_{\rm Fr}(\g_{\bar{0}},\F) $
and
$H: H^2(\g)\to\Hom_{\rm Fr}(\g_{\bar{0}},H^1(\g))$
in
(\ref{sixterm}) are  given by

\begin{equation}\label{D}
D_\psi(g)= \psi(g^{[p]})
  \end{equation}
and
\begin{equation}\label{HH}
H_\phi(g)\cdot h =\varphi(g\wedge (\ad g)^{p-1}(h))-\phi(g^{[p]}\wedge h)
  \end{equation}
where $g\in \g_{\bar{0}}$ and $h\in\g$.
\begin{proof}
The proof is similar to Theorems 2.1, 3.1 and 3.2 in \cite{H}.
\end{proof}

\subsection{Restricted one-dimensional central extension}
The definition of restricted central
extensions was introduced in \cite{YYC}.
Let $(\g,[p])$
be a restricted Lie superalgebra and $M$ be a strongly abelian restricted Lie superalgebra (i.e., $[M, M]=0$ and $M_{\bar{0}}^{[p]}=0$). 
A \emph{restricted extension} of $\g$ by $M$ is a short exact sequence of restricted Lie superalgebras
\begin{align*}
  0 \longrightarrow   M \stackrel{\iota}{\longrightarrow}  \mathfrak{G} \stackrel{\pi}{\longrightarrow}  \g \longrightarrow  0.
\end{align*}
In the above diagram, the maps  $\iota$ and $\pi$ are even restricted Lie superalgebra homomorphisms. If $\iota(M)$ is contained in the center of $\mathfrak{G}$, then the extension is called \emph{central}. Two restricted central extensions
of $\g$ by $M$ are called \emph{equivalent} if they can be included in the usual
commutative diagram. That is, there exists a restricted Lie superalgebra homomorphism $\sigma: \mathfrak{G}_1\to \mathfrak{G}_2$ that ﬁxes $M$ elementwise and $\pi_2 \sigma=\pi_1$.

Let $(\g,[p])$ be a restricted Lie superalgebra and $M$ a restricted $\g$-module. We define a subspace
of $C^2_{res}(\g,M)$ by
$$C^2_{res}(\g,M)^{+}=\{(\phi,\omega)\in C^2_{res}(\g,M)\mid \omega(\g_{\bar{0}}) \subseteq M_{\bar{0}} \}.$$
Since $B_{res}^2(\g,M)_{\bar{0}}=d_{res}^{1}(C^1_{res}(\g,M)_{\bar{0}})\subseteq
C^2_{res}(\g,M)^{+}$, we define
$$Z_{res}^{2}(\g,M)^{+}=\mathrm{ker}\left(d^2_{res\mid C^2_{res}(\g,M)^{+}}\right)$$ and
$$H_{res}^{2}(\g,M)^{+}=Z_{res}^{2}(\g,M)^{+}/B_{res}^2(\g,M)_{\bar{0}}.$$ In the case of $M=\F$, we have
$C^2_{res}(\g)^{+}=C^2_{res}(\g)$
and
$H_{res}^{2}(\g)^{+}_{\bar{0}}=H_{res}^{2}(\g)_{\bar{0}}.$ Moreover, we have the following theorem.
\begin{theorem}\cite[Theorem 3.5.3]{BE}\label{ce}
Let $(\g,[p])$ be a restricted Lie superalgebra. Then $H_{res}^{2}(\g)_{\bar{0}}$ is in one to one correspondence with the equivalence classes of restricted one-dimensional central extensions of $\g$.
\end{theorem}
Now we consider restricted central 
extensions of $\g$ by a one-dimensional space $\F c$. By Theorem \ref{ce}, these restricted one-dimensional central extensions can be determined by restricted 2-cohomology.
In particular,
if $(\phi,\omega) \in C^2_{res}(\g)_{\bar{0}}$ is a
restricted even 2-cocycle, then the corresponding restricted
one-dimensional central extension $\mathfrak{G}=\g\oplus\F c$ has  the bracket and
$[p]$-operator defined for all $g,h\in\g$ and $g_0\in\g_{\bar{0}}$ by
\begin{align}\label{genonedimext}
  \begin{split}
  [g,h]_{\mathfrak{G}}&=[g,h]_{\g} + \phi(g\wedge h)c\\
  g_0^{[p]_{\mathfrak{G}}}&=g_0^{[p]_{\g}} + \omega(g_{0})c
  \end{split}
\end{align}
where $[\cdot,\cdot]_{\g}$ and $[p]_{\g}$ denote the 
bracket of $\g$ and $[p]$-operator of $\g$, respectively.

\section{Restricted Lie superalgebra $\h^{\lambda}_{m,n}$}
Recall the Heisenberg
Lie superalgebra with the even center
$\h_{m,n}$ having a basis 
\[\{x_1,\ldots , x_{2m},x_{2m+1}\mid y_1,\ldots, y_n\}\] 
with non-zero brackets of basis
elements (in particular, $x_{2m+1}$ spans the even center),
\[[x_i,x_{m+i}]=[y_j,y_j]=x_{2m+1}\] 
for $m,n\geq 1$, $1\le i\le m$ and $1\le j\le n$.
The even part of $\h_{m,n}$ is exactly  the Heisenberg Lie algebra $\h_m$ \cite[Section 3]{R} and the restricted structures of $\h_m$, parameterized by
 elements $\lambda\in \F^{2m+1}$,  were studied in \cite[Section 3.1]{EFY}.  
 If 
 $$\lambda=(\lambda_1,\dots, \lambda_{2m+1})\in \F^{2m+1}$$ 
 and  
 $g= \displaystyle\sum^{2m+1}_{i=1} a_i x_i\in (\h_{m,n})_{\bar{0}}=\h_m$,
then
\begin{equation}\label{res1}
  g^{[p]}=\left(\sum^{2m+1}_{i=1}a^{p}_i\lambda_i\right) x_{2m+1}.
\end{equation}
Since $[(\h_{m,n})_{\bar{0}},(\h_{m,n})_{\bar{1}}]=0$, $(\h_{m,n},[p])$ is a restricted Lie superalgebra by Definition \ref{res-Lie-super}.
We denote the
corresponding restricted Lie superalgebra by $\h^{\lambda}_{m,n}$.

\subsection{Restricted cohomology $H^q_{res}(\h_{m,n}^\lambda)$ for $q=1,2$}
The cohomology of the Heisenberg Lie superalgebras $\h_{m,n}$ with trivial
coefficients in characteristic $0$ was studied in \cite[Theorem 4.1]{BaLi} and it is easy to check the results on $H^{1}(\h_{m,n})$ and $H^{2}(\h_{m,n})$
still hold for characteristic $p>2$:

\begin{lemma}\label{Hh11}
  For $p>2$,

  (1) 
  $H^{1}(\h_{m,n})=(\h_{m,n}/\F x_{2m+1})^{\ast};$

  (2) $
    H^{2}(\h_{m,n})=\bigwedge^{2}(\h_{m,n}/\F
                 x_{2m+1})^{\ast}/(\F dx^{2m+1}) 
$
  where $d x^{2m+1}=\sum^{m}_{i=1} x^{i,m+i}-2^{-1}\sum_{i=1}^{n}  y^{i,i}$.
\end{lemma}

Let us now consider the restricted Lie superalgebra
$\mathfrak{h}_{m,n}^{\lambda}$. The
$[p]$-operator formula (\ref{res1})  implies
$\langle (\h^\lambda_{m,n})_{\bar{0}}^{[p]}\rangle=\F x_{2m+1}$. Then, it follows from Equation (\ref{h1}) that 
 $$H^1(\h_{m,n}^\lambda)=H^1_{res}(\h_{m,n}^\lambda),$$
 the
classes of $\{x^1,\ldots , x^{2m} \mid y^1,\ldots, y^n\}$ form a basis and the map
$H^1_{res}(\h_{m,n}^\lambda) \to H^1(\h_{m,n}^\lambda)$ is an isomorphism.
Besides that,  Equations (\ref{HH}) and (\ref{res1}) give $H=0$ in the sequence (\ref{sixterm})
 and the
six-term exact sequence (\ref{sixterm}) reduces
 to the exact
sequence
\begin{align*}
  0\longrightarrow\Hom_{\rm Fr}((\h_{m,n}^\lambda)_{\bar{0}},\F) \longrightarrow H^2_{res}(\h_{m,n}^\lambda) \longrightarrow H^2(\h_{m,n}^\lambda)\longrightarrow 0.
\end{align*}
It follows 
 $H^2_{res}(\h_{m,n}^\lambda)\cong\Hom_{\rm Fr}((\h_{m,n}^\lambda)_{\bar{0}},\F)\oplus
   H^2(\h_{m,n}^\lambda).$ Part (2) of Lemma \ref{Hh11} implies   the following theorem:
\begin{theorem}
  \label{maintheorem}
  For $m,n\ge 1$ and any form $\lambda:(\h_{m,n})_{\bar{0}}\to\F$,
  $\mathrm{sdim}\ H^2_{res}(\h_{m,n}^\lambda)=(2m^2+m+\frac{n^2+n}{2},2mn)$ and a basis consists of the
  classes of the cocycles
  \[\{(x^{i,j},\widetilde{x^{i,j}})
   \ |\ 1\le i< j\le 2m\}
\bigcup \{(x^{i}\wedge y^j,\widetilde{x^{i}\wedge y^j})\ |\ 1\le i\le 2m, 1\le  j\le n\}\]
\[\bigcup \{(y^{i,j},\widetilde{y^{i,j}})\ |\ 1\le i< j\le n\}
  \bigcup \{(y^{i,i},\widetilde{y^{i,i}})\ |\ 1\le i\le n-1\}
    \bigcup \{(0,\overline x^i)\ |\ 1\le i\le 2m+1\}.\]
\end{theorem}

\subsection{Restricted one-dimensional central extensions of
  $\h_{m,n}^\lambda$}
With the Equations (\ref{genonedimext}) together with
Theorem~\ref{maintheorem} we can explicitly describe the restricted
one-dimensional central extensions of $\h_{m,n}^\lambda$. 
Let $g=\sum a_ix_i+\sum b_iy_i$, $h=\sum a_i' x_i+ \sum b'_iy_i$ and $g_{0}=\sum c_ix_i$ denote  arbitrary elements of
$\h_{m,n}^\lambda$ and $(\h_{m,n}^\lambda)_{\bar{0}}$, respectively.

If $1\le i\le 2m+1$ and $\mathfrak{H}_i=\h_{m,n}^\lambda\oplus \F c$ denotes
the one-dimensional restricted central extension of $\h_{m,n}^\lambda$ by a one-dimensional space
$\F c$
determined by the cohomology class of the restricted cocycle
$(0,\overline x^i)$, then Equations (\ref{genonedimext}) gives the (non-zero) bracket   and
$[p]$-operator 
\begin{align*}
  \begin{split}
    [g,h]_{\mathfrak{H}_i} & =[g,h]_{\h_{m,n}^\lambda};\\
    g_0^{[p]_{\mathfrak{H}_i}} & = g^{[p]_{\h_{m,n}^\lambda}}_{0}+ c_i^p c .
  \end{split}
\end{align*}
The central extensions $\mathfrak{H}_i$ form a basis for the $(2m+1)$-dimensional
space of restricted one-dimensional central extensions that split as
ordinary Lie superalgebra extensions (c.f. \cite{EFi2}).
 
Since $p>2$, the $(p-1)$-fold bracket in Equation (\ref{starprop}) is a multiple
of $x_{2m+1}$. If $1\le s< t\le 2m$ and $1\le k\le l\le n$, the maps $x^{s,t}$ and $y^{k,l}$ vanish on
$x_{2m+1}\wedge \h_{m,n}^\lambda$ so that Equation (\ref{starprop}) implies that 
$\widetilde{x^{s,t}}$ and  $\widetilde{y^{k,l}}$ are $p$-semilinear, and hence
$\widetilde{x^{s,t}}=\widetilde{y^{k,l}}=0$. It follows that for any basis element in
Theorem~\ref{maintheorem} of the form
$(x^{s,t},\widetilde{x^{s,t}})=(x^{s,t},0)$ and $(y^{k,l},\widetilde{y^{k,l}})=(y^{k,l},0)$. The corresponding restricted one-dimensional
central extension $\mathfrak{X}_{s,t}$ determined by the cohomology class of the restricted cocycle $(x^{s,t},0)$ has (non-zero) bracket and
$[p]$-operator
\begin{align*}
  \begin{split}
    [g,h]_{\mathfrak{X}_{s,t}} & =[g,h]_{\h_{m,n}^\lambda}+(a_sa'_t-a_ta'_s)c;\\
    g_{0}^{[p]_{\mathfrak{X}_{s,t}}} & = g_{0}^{[p]_{\h_{m,n}^\lambda}}.
  \end{split}
\end{align*}
The corresponding restricted one-dimensional
central extension $\mathfrak{Y}_{k,l}$ determined by the cohomology class of the restricted cocycle $(y^{k,l},0)$ has (non-zero) bracket and
$[p]$-operator
\begin{align*}
  \begin{split}
    [g,h]_{\mathfrak{Y}_{k,l}} & =[g,h]_{\h_{m,n}^\lambda}-(b_kb'_l+b_lb'_k)c;\\
    g_{0}^{[p]_{\mathfrak{Y}_{k,l}}} & = g_{0}^{[p]_{\h_{m,n}^\lambda}}.
  \end{split}
\end{align*}
The central extensions $\mathfrak{X}_{s,t}$ and $\mathfrak{Y}_{k,l}$ form a basis for the $(2m^2-m+\frac{n^2+n}{2}-1)$-dimensional
space of restricted one-dimensional central extensions that do not split as
ordinary Lie superalgebra extensions and the Lie algebras $(\mathfrak{Y}_{k,l})_{\bar{0}}$ cannot be obtained from the central extensions of Heisenberg Lie algebras \cite{EFY}.

\section{Restricted Lie superalgebra $\mathfrak{ba}^{0}_n$}
Recall the Heisenberg Lie superalgebra with the odd center $\mathfrak{ba}_n$ having a basis
\[\{x_1,\ldots , x_{n}\mid y_1,\ldots, y_n,y_{n+1}\}\] with
non-zero brackets of basis elements (in particular, $y_{n+1}$ spans the odd center),
  \[[x_i,y_i]=y_{n+1}\] 
  for $n\geq 1$ and $1\le i\le n$.
Suppose that $[p]$ is a $[p]$-operator on $(\mathfrak{ba}_n)_{\bar{0}}$ such that $(\mathfrak{ba}_n,[p])$ is a restricted Lie superalgebra.
Then $[(\mathfrak{ba}_n)_{\bar{0}}^{[p]},(\mathfrak{ba}_n)_{\bar{1}}]=0$ by Definition \ref{res-Lie-super}. Then $(\mathfrak{ba}_n)_{\bar{0}}^{[p]}=0$. That is, $[p]=0$. It means that the only $[p]$-operator is 0. 
 We denote the
corresponding restricted Lie superalgebra by $\mathfrak{ba}^{0}_{n}$.

\subsection{Restricted cohomology $H_{res}^q(\mathfrak{ba}^{0}_n)$ for $q=1,2$}
The cohomology of the Heisenberg Lie superalgebras $\mathfrak{ba}_{n}$ with trivial
coefficients in characteristic $0$ was studied in \cite[Theorem 4.2]{BaLi} and it is easy to check the results on $H^{1}(\mathfrak{ba}_{n})$ and $H^{2}(\mathfrak{ba}_{n})$
still hold for characteristic $p>2$:

\begin{lemma}\label{H11}
  For $p>2$,

  (1) 
  $H^{1}(\mathfrak{ba}_{n})=(\mathfrak{ba}_{n}/\F y_{n+1})^{\ast};$

  (2) $
    H^{2}(\mathfrak{ba}_{n})=\bigwedge^{2}(\mathfrak{ba}_{n}/\F
                 y_{n+1})^{\ast}/(\F dy^{n+1}) \oplus \delta_{1,n}(x^1\wedge y^{2}),
$\\
  where  $\delta_{1,n}=1$ when
 $n=1$ and $\delta_{1,n}=0$ otherwise, $d y^{n+1}=\sum^{n}_{i=1} x^{i}\wedge  y^{i}$.
\end{lemma}

Let us now consider the restricted Lie superalgebra
$\mathfrak{ba}_n^{0}$. 
 Since $[p]=0$,
it
follows that $H^1(\mathfrak{ba}_{n}^0)=H^1_{res}(\mathfrak{ba}_{n}^0)$ and $H=0$ in the sequence  (\ref{sixterm}).
With a similar argument to $\h_{m,n}^{\lambda}$, we have the following theorem:
\begin{theorem}
  \label{maintheorem1}
  For $n\ge 1$,
  $\mathrm{sdim}\ H^2_{res}(\mathfrak{ba}_{n}^0)=(n^2+n,n^2-1+\delta_{1,n})$ and a basis consists of the
  classes of the cocycles
  \[\{(x^{i,j},\widetilde{x^{i,j}})
   \ |\ 1\le i< j\le n\}
\bigcup \{(y^{i,j},\widetilde{y^{i,j}})\ |\ 1\le i\le j\le n\}
    \bigcup \{(0,\overline x^i)\ |\ 1\le i\le n\}.\]
    \[\bigcup \{(x^{j}\wedge y^i,\widetilde{x^{j}\wedge y^i})\ |\ 1\le i<j\le n\}
 \bigcup \{(x^{i}\wedge y^i,\widetilde{x^{i}\wedge y^i})\ |\ 1\le i\le n-1\}   
    \]
\[\bigcup \{\delta_{1,n}(x^{1}\wedge y^2,\widetilde{x^{1}\wedge y^2}),(x^{1}\wedge y^i,\widetilde{x^{1}\wedge y^i})\ |\ 2\le i\le n\}
 \bigcup \{(x^{i}\wedge y^j,\widetilde{x^{i}\wedge y^j})\ |\ 2\le i<j\le n\}.   
    \]   
\end{theorem}

\subsection{Restricted one-dimensional central extensions of
  $\mathfrak{ba}^{0}_n$}

Let $g=\sum a_ix_i+\sum b_iy_i$, $h=\sum a_i' x_i+ \sum b'_iy_i$ and $g_{0}=\sum c_ix_i$ denote  arbitrary elements of
$\mathfrak{ba}^{0}_n$ and $(\mathfrak{ba}^{0}_n)_{\bar{0}}$, respectively.

If $1\le i\le n$ and $\mathfrak{H}_i=\mathfrak{ba}^{0}_n\oplus \F c$ denotes
the one-dimensional restricted central extension of $\mathfrak{ba}^0_n$ by a one-dimensional space $\F c$
determined by the cohomology class of the restricted cocycle
$(0,\overline x^i)$, then Equations (\ref{genonedimext}) gives the (non-zero) bracket  and
$[p]$-operator 
\begin{align*}
  \begin{split}
    [g,h]_{\mathfrak{H}_i} & =[g,h]_{\mathfrak{ba}^{0}_n};\\
    g_0^{[p]_{\mathfrak{H}_i}} & = g^{[p]_{\mathfrak{ba}^{0}_n}}_{0}+ c_i^p c.
  \end{split}
\end{align*}
The central extensions $\mathfrak{H}_i$ form a basis for the $n$-dimensional
space of restricted one-dimensional central extensions that split as
ordinary Lie superalgebra extensions (c.f. \cite{EFi2}).
 
Since $p>2$, the $(p-1)$-fold bracket in Equation (\ref{starprop}) is a multiple
of $y_{n+1}$. If $1\le s\le t\le n$, the maps $x^{s,t}$ and $y^{s,t}$ vanish on
$y_{n+1}\wedge \mathfrak{ba}^0_{n}$ so that Equation (\ref{starprop}) implies that 
$\widetilde{x^{s,t}}$ and  $\widetilde{y^{s,t}}$ are $p$-semilinear, and hence
$\widetilde{x^{s,t}}=\widetilde{y^{s,t}}=0$. It follows that for any basis element in
Theorem~\ref{maintheorem1} of the form
$(x^{s,t},\widetilde{x^{s,t}})=(x^{s,t},0)$ and $(y^{s,t},\widetilde{y^{s,t}})=(y^{s,t},0)$. The corresponding restricted one-dimensional
central extension $\mathfrak{X}_{s,t}$ determined by the cohomology class of the restricted cocycle $(x^{s,t},0)$ has (non-zero) bracket and
$[p]$-operator
\begin{align*}
  \begin{split}
    [g,h]_{\mathfrak{X}_{s,t}} & =[g,h]_{\mathfrak{ba}_{n}^0}+(a_sa'_t-a_ta'_s) c;\\
    g_{0}^{[p]_{\mathfrak{X}_{s,t}}} & = g_{0}^{[p]_{\mathfrak{ba}_{n}^0}}.
  \end{split}
\end{align*}
The corresponding restricted one-dimensional
central extension $\mathfrak{Y}_{s,t}$ determined by the cohomology class of the restricted cocycle $(y^{s,t},0)$ has (non-zero) bracket and
$[p]$-operator
\begin{align*}
  \begin{split}
    [g,h]_{\mathfrak{Y}_{s,t}} & =[g,h]_{\mathfrak{ba}_{n}^0}-(b_sb'_t+b_tb'_s)c;\\
    g_{0}^{[p]_{\mathfrak{Y}_{s,t}}} & = g_{0}^{[p]_{\mathfrak{ba}_{n}^0}}.
  \end{split}
\end{align*}
The central extensions $\mathfrak{X}_{s,t}$ and $\mathfrak{Y}_{s,t}$ form a basis for the $n^2$-dimensional
space of restricted one-dimensional central extensions that do not split as
ordinary Lie superalgebra extensions.
\\

\small\noindent \textbf{Acknowledgment}\\
The author thanks the referee for the suggestions which improved the presentation.

\end{document}